\documentclass[11pt]{amsart}
\usepackage[T1]{fontenc}
\usepackage[latin1]{inputenc}

\makeatletter

\providecommand{\LyX}{L\kern-.1667em\lower.25em\hbox{Y}\kern-.125emX\@}

 \theoremstyle{plain}    
 \newtheorem{thm}{Theorem} 
 \theoremstyle{definition}
 \newtheorem*{defn*}{Definition}

\theoremstyle{definition}
\newtheorem*{notations*}{Notations}
\newtheorem*{Remarks}{Remarks}

\newcommand{\ZZ}{\mathbb {Z}}
\newcommand{\TT}{\mathbb {T}}
\newcommand{\RR}{\mathbb {R}}
\newcommand{\CC}{\mathbb {C}}

\makeatother

\begin{document}

\title{An ``Analytic'' Version of Menshov's Representation Theorem}

\author{Gady Kozma, Alexander Olevski\v\i}

\address{School of Mathematical Sciences, Tel Aviv University, Ramat-Aviv 69978, Israel}

\email{\{gady,olevskii\}@math.tau.ac.il}

\begin{abstract}
Every measurable function \( f \) on the circle \( \TT  \) can be represented
as a sum of harmonics with positive spectrum, 
\[
f(t)=\sum _{k>0}c_{k}e^{ikt}\]
converging in measure. For convergence almost everywhere this is not true. We
discuss several other sets \( \Lambda \subset \ZZ  \) for which one might get
a Menshov type representation converging almost everywhere or in measure.
\end{abstract}
\maketitle

\section{Menshov's Theorem}

We denote by \( L^{0}(\TT ) \) the space of measurable functions \( f:\TT \rightarrow \CC  \).
D. E. Menshov proved (1940) that every \( f\in L^{0}(\TT ) \) can be expanded
in a trigonometric series
\begin{equation}
\label{f_is_sum}
f(t)=\sum _{k\in \ZZ }c_{k}e^{ikt}
\end{equation}
converging almost everywhere (a.e.), see {[}1{]}. This famous result simulated
many further investigations. We mention here only papers {[}2-6{]}. The reader
may find a comprehensive bibliography in the survey {[}7{]}.

The expansion (\ref{f_is_sum}) in general is in no sense the Fourier expansion
of \( f \). The coefficients \( \{c_{k}\} \) are defined by a special construction,
involving some free parameters; this makes the expansion far from unique. In
particular one may use high frequency harmonics only. Even infinitely many frequencies
can be removed.

\section{Lacunary Menshov Spectra}

We introduce the following

\begin{defn*}
A set \( \Lambda  \) of integers, \( \Lambda =\{\lambda (n)\, ;\, \ldots <\lambda (-1)<\lambda (0)<\lambda (1)<\ldots \} \)
is a \emph{Menshov spectrum} if every \( f\in L^{0}(\TT ) \) admits a representation
\[
f(t)=\sum _{k\in \Lambda }c_{k}e^{ikt}\equiv \lim _{N\rightarrow \infty }\sum _{k\in \Lambda ,|k|\leq N}c_{k}e^{ikt}\]
converging a.e.
\end{defn*}
It is well known that a Menshov spectrum may have density zero. The following
version of the result is due to Arutyunyan {[}4{]}: if \( \Lambda  \) is a
symmetric set with respect to zero containing arbitrarily long segments of integers
then it is a Menshov spectrum.

Our first result shows that Menshov spectra may even be lacunary, that is satisfy
the condition
\[
\lambda (n+1)-\lambda (n)\rightarrow \infty \quad (|n|\rightarrow \infty )\]
 and the sizes of the gaps may grow quite fast. Certainly, it is impossible
to achieve Hadamarian lacunarity in a Menshov spectrum, but one might get arbitrarily
close.

\begin{thm}
Given a positive sequence \( \epsilon (n)=o(1) \) as \( n\rightarrow \infty  \),
one can define a Menshov spectrum \( \Lambda  \) such that
\[
\frac{\lambda (-n-1)}{\lambda (-n)},\frac{\lambda (n+1)}{\lambda (n)}>1+\epsilon (n)\quad n=1,2,\ldots \]
 
\end{thm}
For specific \( \Lambda  \), even if the gaps grow slowly, it might be difficult
to identify Menshov spectra; some arithmetic obstacles may appear. The following
result concerns perturbation of squares:

\begin{thm}
For any sequence \( w(k)\rightarrow \infty  \) as \( k\rightarrow \infty  \)
one can construct a Menshov spectrum \( \Lambda  \),
\[
\Lambda =\left\{ \pm k^{2}+o(w(|k|))\right\} \quad .\]

\end{thm}
This result is sharp, the reminder cannot be replaced by \( O(1) \).

\section{The Analytic Case}

Now we consider \( \Lambda =\ZZ ^{+} \). Our first observation is that \emph{it
is not a Menshov spectrum}. Indeed, having an expansion
\[
f(t)=\sum _{k>0}c_{k}e^{ikt}\quad \mathrm{a}.\mathrm{e}.\; \mathrm{on}\; \TT ,\]
consider the corresponding analytic function F on the disc \( |z|<1 \)
\[
F(z)=\sum _{k>0}c_{k}z^{k}\quad .\]
 For almost every \( t\in \TT  \), \( F(z)\rightarrow f(t) \) when \( z \)
approaches the point \( e^{it} \) non-tangentially. So, according to Lusin-Privalov
uniqeness theorem, see for example {[}8{]}, \( f \) can not vanish on a set
of positive measure unless it is identically zero. It turns out, however, that
replacing pointwise convergnce by convergence in the \( L^{0} \) metric (that
is, convergence in measure) we may again get an unrestricted representation
theorem

\begin{thm}
\label{Theorem_L0}Every \( f\in L^{0}(\TT ) \) can be represented as a sum
\begin{equation}
\label{f_is_sum_measure}
f(t)=\sum _{k>0}c_{k}e^{ikt}
\end{equation}
converging in measure.
\end{thm}
Actually, we have the result in the following stronger form:

\begin{thm}
Given \( f\in L^{0}(\TT ) \) one can define coefficients \( \{c_{k}\} \) and
compacts \( K_{n}\subset \TT  \), \( \mathbf{m}\left( \TT \setminus K_{n}\right) \rightarrow 0 \),
such that the expansion \emph{(\ref{f_is_sum_measure})} holds in \( L^{2}(K_{n}) \)
for all \( n \).
\end{thm}
In addition one can require sharp conditions on the decrease of the coefficients
\( c_{k} \) (certainly, beyond \( l_{2} \)), for example \( \{c_{k}\}\in l_{2+\epsilon } \)
for any \( \epsilon >0 \).

\begin{Remarks}

\begin{enumerate}
\item In the theorem above one can replace convergence in \( L^{2} \) of \( K_{n} \)
by \( L^{p} \)-convergence for any \( p<\infty  \) (but not for \( p=\infty  \)).
Such kind of convergence for any \( 1\leq p<\infty  \) was considered by Talalyan
{[}3{]} in his representation results with respect to complete orthonormal systems
and bases. He used the term ``asymptotic convergence in \( L^{p} \)'' for
it.
\item Theorem \ref{Theorem_L0} admits also functions \( f:\TT \rightarrow \RR \cup \{\pm \infty \} \).
This improves the Menshov representation theorem for such functions (see {[}2{]})
by replacing the whole spectrum \( \ZZ  \) by \( \ZZ ^{+} \).
\item It might be interesting to compare theorems 3 and 4 with the ``radial representation
theorem'' of Kahane and Katznelson {[}9{]}, where \( f \) is obtained as the
radial limit a.e. of an analytic function.
\end{enumerate}
The ``lacunary'' versions of the theorems of this section are also true. For
example:\end{Remarks}

\begin{thm}
Given \( 0<\epsilon (n)=o(1) \) one can define \( \Lambda \subset \ZZ ^{+} \),
\( \frac{\lambda (n+1)}{\lambda (n)}>1+\epsilon (n) \) such that every \( f\in L^{0}(\TT ) \)
admits a representation
\begin{equation}
\label{f_is_sum_lambda}
f(t)\mathop {=}^{L^{0}}\sum _{k\in \Lambda }c_{k}e^{ikt}
\end{equation}
 
\end{thm}
{}

\begin{thm}
For any sequence \( w(k)\rightarrow \infty  \) as \( k\rightarrow \infty  \)
one can construct a set \( \Lambda \subset \ZZ ^{+} \), 
\[
\Lambda =\left\{ k^{2}+o(w(|k|))\right\} \quad .\]
such that every \( f\in L^{0}(\TT ) \) admits a representation (\ref{f_is_sum_lambda})
\end{thm}

Again, the result is sharp: it is not true for a bounded perturbation of the
squares. We conclude with the following result in the spirit of Arutyunyan {[}4{]}

\begin{thm}
If a set \( \Lambda \subset \ZZ ^{+} \) contains arbitrarily long segments
then \( \forall f\in L^{0} \) can be represented as in \emph{(\ref{f_is_sum_lambda})}.
\end{thm}

\end{document}